%
%
 
 \magnification=1200
 \tolerance=999
\overfullrule=0pt    \parskip=13pt \parindent=15pt

 
\catcode`@=11 
 
\font\sixrm=cmr7 scaled 833
\font\sixbf=cmbx7 scaled 833
 
\font\eighti=cmmi10 scaled 833
\font\sixi=cmmi7 scaled 833
 
\font\eightsy=cmsy10 scaled 833
\font\sixsy=cmsy7 scaled 833

\font\eightbf=cmbx10 scaled 833
 
\font\eighttt=cmtt10 scaled 833
 
\font\eightsl=cmsl10 scaled 833
 
\font\eightit=cmti10 scaled 833
 
\font\fiverm=cmr5 scaled 833
\font\fivei=cmmi5 scaled 833
\font\fivesy=cmsy5 scaled 833
 
\font\tenux=cmex10 scaled 833
 
\def\eightpoint{\def\rm{\fam0\fsmall}%
  \textfont0=\fsmall \scriptfont0=\sixrm \scriptscriptfont0=\fiverm
  \textfont1=\eighti \scriptfont1=\sixi \scriptscriptfont1=\fivei
  \textfont2=\eightsy \scriptfont2=\sixsy \scriptscriptfont2=\fivesy
  \textfont3=\tenux \scriptfont3=\tenux \scriptscriptfont3=\tenux
  \def\it{\fam\itfam\eightit}%
  \textfont\itfam=\eightit
  \def\sl{\fam\slfam\eightsl}%
  \textfont\slfam=\eightsl
  \def\bf{\fam\bffam\eightbf}%
  \textfont\bffam=\eightbf \scriptfont\bffam=\sixbf
   \scriptscriptfont\bffam=\fivebf
  \def\tt{\fam\ttfam\eighttt}%
  \textfont\ttfam=\eighttt
  \normalbaselineskip=9pt
  \def\MF{{\manual opqr}\-{\manual stuq}}%
  \let\sc=\sixrm
  \let\big=\eightbig
  \setbox\strutbox=\hbox{\vrule height7pt depth2pt width\z@}%
  \normalbaselines\rm}

\def\footnote#1{ \parindent=15pt
                \edef\@sf{\spacefactor\the\spacefactor}#1\@sf
      \insert\footins\bgroup\eightpoint
      \interlinepenalty100 \let\par=\endgraf
        \leftskip=\z@skip \rightskip=\z@skip
        \splittopskip=10pt plus 1pt minus 1pt \floatingpenalty=20000
        \smallskip\item{#1}\bgroup\strut\aftergroup\@foot\let\next}
\skip\footins=12pt plus 2pt minus 4pt 
\dimen\footins=30pc 

\catcode`@=12
 
\def \nheadline{\hfil \tenrm \the\pageno \hfil}
\footline={}
 
\newif \iftitlepage \titlepagetrue
\headline=
  {\iftitlepage \hfil \global \titlepagefalse \else \nheadline \fi}

\newcount\footnotenum \footnotenum=0
\def\fnote#1{\advance \footnotenum by 1%
           \showthe\footnotenum
\footnote{$^{\the\footnotenum}$}{ \fsmall{#1} }
            }

\font\fsmall=cmr10 scaled 833
\font\ttt=cmtt10 scaled 833
\font\fgrand=cmbx12
\font\blackb=msbm10 scaled 833

\def\dfrac#1#2{{\displaystyle{#1\over#2}}}

\def\R{{\hbox{\blackb R}}}
 

{
{
\baselineskip=14pt
\noindent
\centerline{\fsmall Stieltjes meeting, Delft, 31 Oct.-4 Nov. 1994}
\centerline{\fgrand Problem}
\centerline{\fsmall
 for Mathematicaists, Reducectionists, Macsymalists, and other Mapleators}
\centerline
{\fgrand Painlev\'e equations for semi-classical recurrence coefficients.}
}
 
\vskip  2mm
 
\centerline{Alphonse P. Magnus}
 
\centerline
{\fsmall Institut Math\'ematique,  Universit\'e Catholique de Louvain,
 Chemin du Cyclotron 2,
 B-1348 Louvain-la-Neuve,
 Belgium.}
 
\centerline{\fsmall e-mail:  {\ttt magnus@anma.ucl.ac.be}}

\bigskip
 
}  
 
\displayindent=5pt
 
{\bf Background:\/}
orthonormal polynomials $p_n(x)=\gamma_n x^n+\cdots\ \ $:
$$ a_{n+1}p_{n+1}(x)=(x-b_n)p_n(x)-a_n p_{n-1}(x) \eqno(1)$$
related to a semi-classical weight\ $w:  \qquad \qquad
   \dfrac{w'(x)}{w(x)}=\dfrac{2V(x)}{W(x)}$  \hfill(2)\break
where $V$ and $W$ are polynomials, satisfy differential relations
$$Wp_n' = (\Omega_n-V)p_n -a_n \Theta_n p_{n-1} \eqno(3)$$
where $\Omega_n$ and $\Theta_n$ are polynomials of degrees $\le m-1$
and $m-2$ (if degrees $V$ and $W \le m-1$ and $m$), with
\qquad
$ \Omega_{n+1}(x)=(x-b_n)\Theta_n(x)-\Omega_n(x), \hfill (4)$
$$ (x-b_n)(\Omega_{n+1}(x)-\Omega_n(x))=W(x)+a_{n+1}^2\Theta_{n+1}(x)-
                                             a_n^2\Theta_n(x)$$
representing nonlinear recurrence relations for the $a_n$'s and
$b_n$'s [Laguerre, $19^{\rm th}$ century].
 
{\bf Statement:\/} {\sl if, in addition to $(2)$, the weight $w$
depends on a parameter $t$ so that
$$ \dfrac{\partial w(x,t)/\partial t}{w(x,t)} =
 \dfrac{T(x,t)}{W(x,t)}, \eqno(5)$$
where $T$ and $W$ are polynomials in $x$ and $t$ $($ the same $W$,
which will now be written $W(x,t))$, each recurrence coefficient
$a_n^2(t)$ and $b_n(t)$ satisfy a nonlinear differential equation having
Painlev\'e property $($movable singular points are only poles$)$ .
} [many people, $20^{\rm th}$ century].
 
I have an (incomplete) elementary proof of this for {\bf generalized
Jacobi weights\/}
$$\left.
  \eqalign{
     w(x,t) &= 0, \qquad x<x_1(t) {\rm\ or\ } x>x_m(t) \cr
  &= C_j \prod_{k=1}^m |x-x_k(t)|^{\alpha_k},
  \quad x_j(t)<x<x_{j+1}(t), \;\;j=1,\ldots,m-1.\cr
          }   \right\}\eqno(6)
$$
where therefore $W(x,t)=\prod_1^m(x-x_k(t))$ and (2) and (5)
become
$$ \dfrac{w'(x,t)}{w(x,t)}=\sum_{k=1}^m \dfrac{\alpha_k}
                                            {x-x_k(t)}
  {\rm\ and\ }
   \dfrac{\partial w(x,t)/\partial t}{w(x,t)}=
  -\sum_{k=1}^m \dfrac{\alpha_k \dot x_k(t)} {x-x_k(t)}
\eqno(7)
$$
(N.B. so, the $C_j$'s and the $\alpha_j$'s do {\it not\/} depend
on $t$; remark also that $\alpha_k=2V(x_k)/W'(x_k)$,
$\Omega_n(x)=(n+(\sum_1^m \alpha_k)/2)x^{m-1}+\cdots,\
\Theta_n(x)=(2n+1+\sum_1^m\alpha_k)x^{m-2}+\cdots $).
A differential system involving $a_n$ and $b_n$ is then
$$ \dfrac{\dot a_n}{a_n} = \dfrac12 \sum_{k=1}^m
    \dfrac{(\Theta_n(x_k)-\Theta_{n-1}(x_k))\dot x_k}{W'(x_k)},
  \eqno(8) $$
$$ \dot b_n       =          \sum_{k=1}^m
    \dfrac{(\Omega_{n+1}(x_k)-\Omega_n(x_k))\dot x_k}{W'(x_k)}
         = \sum_{k=1}^m
    \dfrac{((x_k-b_n)\Theta_n(x_k)-2\Omega_n(x_k))\dot x_k}{W'(x_k)},
  \eqno(9) $$
$$ \dfrac{\dot\gamma_n}{\gamma_n}=-\dfrac12
  \sum_{k=1}^m \dot x_k \Theta_n(x_k)/W'(x_k)
 \eqno(10) $$
$$ \dfrac d{d t} \dfrac{\Theta_n(x_j)}{W'(x_j)} =
 2\dfrac{\dot\gamma_n}{\gamma_n}  \dfrac{\Theta_n(x_j)}{W'(x_j)} -2
\sum_{k\ne j} \dfrac{\dot x_j -\dot x_k}{x_j-x_k}
    \dfrac{\Theta_n(x_k)\Omega_n(x_j)-\Theta_n(x_j)\Omega_n(x_k)}
      {W'(x_j)W'(x_k)},\eqno(11)$$
$$\hskip-12pt
   \dfrac d{d t} \dfrac{\Theta_{n-1}(x_j)}{W'(x_j)} =
 -2\dfrac{\dot\gamma_{n-1}}{\gamma_{n-1}}  \dfrac{\Theta_{n-1}(x_j)}{W'(x_j)}
  +2
\sum_{k\ne j} \dfrac{\dot x_j -\dot x_k}{x_j-x_k}
 \dfrac{\Theta_{n-1}(x_k)\Omega_n(x_j)-\Theta_{n-1}(x_j)\Omega_n(x_k)}
      {W'(x_j)W'(x_k)}, \eqno(12)$$
$$\dfrac d{d t} \dfrac{\Omega_n(x_j)}{W'(x_j)} =
a_n^2 \sum_{k\ne j} \dfrac{\dot x_j -\dot x_k}{x_j-x_k}
 \dfrac{\Theta_n(x_j)\Theta_{n-1}(x_k)-\Theta_n(x_k)\Theta_{n-1}(x_j)}
      {W'(x_j)W'(x_k)},\eqno(13)$$
 
{\fgrand Problem:\/} eliminate the $\Theta$'s and the $\Omega$'s in
$(8-13)$ so to exhibit the scalar differential equations for $a_n(t)$ and
$b_n(t)$, at least in the $m=3$ case ($W$ of degree 3 in $x$ and 1 in
$t$), and in the even $m=5$ case ($w(x,t)=|x|^\alpha|x^2-x_1^2|^\beta
|x^2-x_2^2|^\gamma, b_n=0$), see [1] for confluent cases and [2] for the
even case.
 
{\bf Proof of } $(8-13)$
 
{\eightpoint
 
The proof is valid only when $\alpha_k >0$, so that
$$ q_n(x)=\int_{\R} \dfrac{w(u)p_n(u)}{x-u} du \eqno(14)$$
is finite at the $x_k$'s. Then, $\dfrac{d}{dt}\int_\R w(u,t)p_n(u,t)
p_{n-i}(u,t)du=0,\quad i=0,1,\ldots$:
$$\sum_{k=1}^m \alpha_k \dot x_k q_n(x_k)p_{n-i}(x_k)
+\int_\R w(u,t){\dfrac{\partial p_n(u,t)}{\partial t}}
p_{n-i}(u,t)du +{\dfrac{\dot\gamma_n}{\gamma_n}}\delta_{i,0}=0,$$
yielding the coefficients of the orthogonal expansion of
$$\displaylines{
 \hskip-12pt
 \dfrac{\partial p_n(x,t)}{\partial t} =
  -\dfrac{\dot\gamma_n}{\gamma_n} p_n(x,t)
  -\sum_{k=1}^m\alpha_k \dot x_k q_n(x_k,t)
  \sum_{i=0}^n p_{n-i}(x_k,t)p_{n-i}(x,t) \hfill(15)\cr
 \hskip 13pt =-\dfrac{\dot\gamma_n}{\gamma_n} p_n(x,t)
  -\sum_{k=1}^m\alpha_k \dot x_k q_n(x_k,t)
 \left[a_n(t)
 \dfrac{p_n(x,t)p_{n-1}(x_k,t)-p_n(x_k,t)p_{n-1}(x,t)}
       {x-x_k} + p_n(x_k,t)p_n(x,t)
    \right] \hfill(16)\cr}      $$
(Christoffel-Darboux). From the coefficient of $p_n$:
$$ \dfrac{\dot\gamma_n}{\gamma_n}=-\dfrac12
  \sum_{k=1}^m \alpha_k \dot x_k q_n(x_k)p_n(x_k)
$$
Now, from (14), (2), (3), (4), and (1)
$$Wq'_n=(\Omega_n+V)q_n-a_n\Theta_n q_{n-1},
Wp'_{n-1}=a_n\Theta_{n-1} p_n -(\Omega_n+V)p_{n-1},
Wq'_{n-1}=a_n\Theta_{n-1} q_n -(\Omega_n-V)q_{n-1},
$$
so that, using $p_nq_{n-1}-p_{n-1}q_n=1/a_n$,
$$\Theta_n(x_k)=2V(x_k)p_n(x_k)q_n(x_k)=\alpha_k W'(x_k)p_n(x_k)q_n(x_k),
\Omega_n(x_k)-V(x_k)=a_n\alpha_k W'(x_k)q_n(x_k)p_{n-1}(x_k),  \eqno(17)
$$
giving (10),
(8) follows from $\gamma_{n-1}=a_n\gamma_n$: $\dot a_n/a_n=
            \dot\gamma_{n-1}/\gamma_{n-1}-\dot\gamma_n/\gamma_n$.
Now, (16) becomes
$$ \dfrac{\partial p_n(x,t)}{\partial t} =
   \dfrac{\dot\gamma_n}{\gamma_n} p_n(x,t)-
 {\displaystyle \sum_{k=1}^m \dot x_k\;\dfrac{(\Omega_n(x_k)-V(x_k))p_n(x)
      -a_n\Theta_n(x_k)p_{n-1}(x)}{W'(x_k)(x-x_k)}}.\eqno(18)$$
For (9), consider
$$\dfrac{\partial p_n(x,t)}{\partial t} =
  \dfrac{\partial [ \gamma_n(t)x^n-\gamma_n(t)(b_0(t)+\cdots
                   +b_{n-1}(t))x^{n-1}+\cdots ]}{\partial t} =
 \dfrac{\dot\gamma_n}{\gamma_n}p_n -\dfrac{\dot b_0+\cdots
        +\dot b_{n-1}}{a_n} p_{n-1}+\cdots $$
and the coefficient of $p_{n-1}$ in (15):
$$-\dfrac{\dot b_0+\cdots+b_{n-1}}{a_n} =
  -\sum_1^m \alpha_k \dot x_k q_n(x_k)p_{n-1}(x_k) =
 -\dfrac1{a_n} \sum_1^m \dot x_k \dfrac{\Omega_n(x_k)-V(x_k)}{W'(x_k)}\;.
$$
 
For (11-13), we need the $t-$derivatives of $p_n(x_j(t),t)$,
                                       $q_n(x_j(t),t)$, etc.
 
$dp_n(x_j(t),t)/dt = (\partial p_n(x,t)/\partial t)(x_j) +
                 \dot x_j p_n'(x_j)$
 
$(3)\Rightarrow\;
 p_n'(x)=\dfrac{\Omega_n(x)-V(x)}{W(x)}p_n(x)-a_n\dfrac{\Theta_n(x)}{W(x)}
  p_{n-1}(x)=
 {\displaystyle \sum_{k=1}^m \dfrac{(\Omega_n(x_k)-V(x_k))p_n(x)
      -a_n\Theta_n(x_k)p_{n-1}(x)}{W'(x_k)(x-x_k)}}$,
and use (18):
$$dp_n(x_j(t),t)/dt =
   \dfrac{\dot\gamma_n}{\gamma_n} p_n(x_j,t) +
     \sum_{\scriptstyle k=1\atop \scriptstyle k\ne j}^m
      (\dot x_j-\dot x_k)\;\dfrac{(\Omega_n(x_k)-V(x_k))p_n(x_j)
      -a_n\Theta_n(x_k)p_{n-1}(x_j)}{W'(x_k)(x_j-x_k)}
    ,\eqno(19)$$
 
$$\eqalign{ \hskip -52pt
  \dfrac{\partial q_n(x,t)}{\partial t} &=
  \dfrac{\partial \int_\R \dfrac{w(u,t)p_n(u,t)}{x-u}du}{\partial t} \cr
 = - &\sum_1^m \alpha_k \dot x_k \int_\R \dfrac{w(u)p_n(u)}{(u-x_k)(x-u)}
   du + \int_\R w(u)\dfrac{\dot\gamma_n p_n(u)}{\gamma_n (x-u)}du
   -\sum_1^m  \dot x_k
  \int_\R w(u)\dfrac{[\Omega_n(x_k)-V(x_k)]p_n(u)-a_n\Theta_n(x_k)p_{n-1}(u)}
                    {W'(x_k)(u-x_k)(x-u)}
   du \cr
 &= \dfrac{\dot\gamma_n}{\gamma_n}q_n(x) -
    \sum_1^m  \dot x_k
  \dfrac{[\Omega_n(x_k)+V(x_k)](q_n(x)-q_n(x_k))
  -a_n\Theta_n(x_k)(q_{n-1}(x)-q_{n-1}(x_k))}{W'(x_k)(x-x_k)}\cr
   }
$$
$$dq_n(x_j(t),t)/dt =
   \dfrac{\dot\gamma_n}{\gamma_n} q_n(x_j,t) +
     \sum_{\scriptstyle k=1\atop \scriptstyle k\ne j}^m
      (\dot x_j-\dot x_k)\;\dfrac{(\Omega_n(x_k)+V(x_k))q_n(x_j)
      -a_n\Theta_n(x_k)q_{n-1}(x_j)}{W'(x_k)(x_j-x_k)}
    ,$$
and we get (11) as $\alpha_j d(p_n(x_j)q_n(x_j))/dt$ (from (17)).
For (12): use (11) with $n-1$ instead of $n$, use $\Omega_{n-1}(x)=
(x-b_{n-1}\Theta_{n-1}(x)-\Omega_n(x)$ from (4), and (12) follows from
(10) and also $\sum_1^m \Theta_{n-1}(x_k)/W'(x_k)=$ coefficient of
$x^{m-1}$ of $\Theta_{n-1}=0$). Finally, we need
$d p_{n-1}(x_j(t),t)/dt$ for (13): take (19) with $n-1$, replace
$-a_{n-1}p_{n-2}(x_j)$ by $a_np_n(x_j)-(x_j-b_{n-1})p_{n-1}(x_j)$:
$$dp_{n-1}(x_j(t),t)/dt =
   -\dfrac{\dot\gamma_{n-1}}{\gamma_{n-1}} p_{n-1}(x_j,t) -
     \sum_{\scriptstyle k=1\atop \scriptstyle k\ne j}^m
      (\dot x_j-\dot x_k)\;\dfrac{(\Omega_n(x_k)+V(x_k))p_{n-1}(x_j)
      -a_n\Theta_{n-1}(x_k)p_n(x_j)}{W'(x_k)(x_j-x_k)}
    .$$

} 
 
{\fgrand Exercise:\/} {\sl show that the moment
$\mu_n=\int_\R w(u) (u-x_1)^n du$ of $(6)$ enters the solution of the
following differential
system as} $\mu_n=\nu_{n,1}$,   \hfill\break
{\sl where\/} $\nu_{n,j}=\int_\R w(u)\prod_1^m(u-x_k)^{\beta_k}
     (u-x_j)^{-1}du$:
 
$$ \dot\nu_{n,j} = \sum_{\scriptstyle k=1\atop \scriptstyle k\ne j}^m
       (\dot x_j-\dot x_k)(\alpha_k+\beta_k)
       \dfrac{\nu_{n,j}-\nu_{n,k}}{x_j-x_k}\,,\qquad
       j=1,\ldots,m $$
{\sl with\/} $\beta_1=n+1,\beta_2=\ldots=\beta_m=1.$
 
This shows that the moments of these generalized Jacobi weights are
regular functions of $t$ in any region where the $x_k(t)$'s are
distinct. As $a_n^2(t)$ and $b_n(t)$ are ratios of determinants of
moments, the movable singular points of their equations can only be
poles (the zeros of the denominator determinants).

\item{[1]} A.P.~MAGNUS, Painlev\'e-type differential equations for
        the recurrence  coefficients of semi-classical orthogonal
         polynomials, to appear in {\sl J.\ Comp.\ Appl.\ Math.}

\item{[2]} A.P.~MAGNUS,
 Asymptotics for the simplest generalized Jacobi polynomials recurrence
  coefficients from Freud's equations: numerical explorations. Preprint.
 
\bye